\documentclass{article}

\PassOptionsToPackage{numbers, compress}{natbib}

\usepackage[final]{neurips_2022}

\usepackage{amsmath}
\usepackage{amsthm}
\usepackage{amssymb}
\usepackage{color}
\usepackage{graphicx}
\usepackage{subcaption}
\usepackage{booktabs}

\usepackage[utf8]{inputenc} 
\usepackage[T1]{fontenc}    
\usepackage{hyperref}       
\usepackage{url}            
\usepackage{booktabs}       
\usepackage{amsfonts}       
\usepackage{nicefrac}       
\usepackage{microtype}      
\usepackage{xcolor}         

\title{Neuro-Symbolic Partial Differential Equation Solver}

\author{
  Pouria A. Mistani$^\dagger$\thanks{corresponding author.} \\
  Nvidia \\
  \texttt{pmistani@nvidia.com} \\
  \And
  Samira Pakravan \thanks{These authors contributed equally to this work.}\\
  UCSB \\
  \texttt{spakravan@ucsb.edu} \\
  \And
  Rajesh Ilango \\
  Nvidia \\
  \texttt{rilango@nvidia.com} \\
  \AND
  Sanjay Choudhry \\
  Nvidia \\
  \texttt{schoudhry@nvidia.com}
  \And  
  Frederic Gibou \\
  UCSB \\
  \texttt{fgibou@ucsb.edu} \\
}

\begin{document}

\maketitle

\vspace{-0.5cm}
\begin{figure}[h]
\centering
\includegraphics[width=0.99\linewidth]{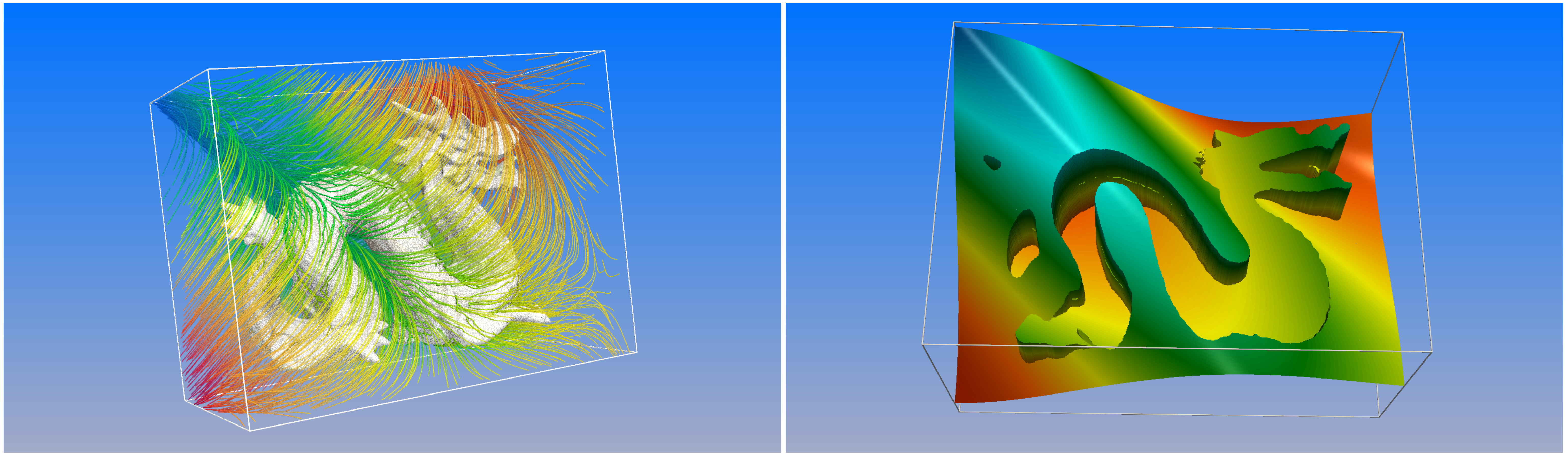}
\caption{The proposed Neural Bootstrapping Method applied to electrostatics. Left: streamlines. Right: warped cross-section showing the physically correct jump in the field. The neural approach readily enables a $1024^3$ effective resolution on a single NVIDIA A6000 GPU. Once trained, it takes sub-milliseconds for the network to evaluate such a simulation \cite{tiny-cuda-nn}.}
\label{fig:dragon}
\end{figure}

\begin{abstract}
 
  We present a highly scalable strategy for developing mesh-free neuro-symbolic partial differential equation solvers from existing numerical discretizations found in scientific computing. This strategy is unique in that it can be used to efficiently train neural network surrogate models for the solution functions and the differential operators, while retaining the accuracy and convergence properties of state-of-the-art numerical solvers. This neural bootstrapping method is based on minimizing residuals of discretized differential systems on a set of random collocation points with respect to the trainable parameters of the neural network, achieving unprecedented resolution and optimal scaling for solving physical and biological systems.
\end{abstract}


\section{Introduction}
Most modern physical, biological and engineering systems are described by partial differential equations on irregular, often moving, boundaries. The difficulties in solving those problems stem from approximating the equations while respecting the physically correct discontinuous nature of the solution across the boundaries. Smoothing strategies are straightforward to design, but unfortunately introduce unphysical characteristics in the solution, which lead to systemic errors.

Since the 1990s, artificial neural networks have been used for solving differential equations using two general strategies: (i) mapping the algebraic operations of the discretized systems onto specialized neural network architectures and minimizing the network energy, or (ii) treating the whole neural network as the basic approximation unit, with parameters trained to minimize a specialized error function that includes the differential equation itself as well as its boundary and initial conditions. 

In the first category, neurons output the discretized solution values over a set of grid points, and minimizing the network energy drives the neuronal values towards the solution at the mesh points. The neural network energy is the residual of the finite discretization, summed over all neurons \cite{lee1990neural}. A strong feature is the preservation of the finite discretization convergence; however, the computational cost grows with increasing resolution and dimensionality. Early examples include \cite{gobovic1993design,chua1988cellular,chua1988cellularA}.

The second strategy, proposed by Lagaris \textit{et al.} \cite{lagaris1998artificial}, relies on the function approximation capabilities of the neural networks. Encoding the solution everywhere in the domain within a neural network offers a mesh-free, compact, and memory efficient surrogate model for the solution function, which can be used in subsequent inference tasks. This method has recently re-emerged as the physics-informed neural networks (PINNs) \cite{RAISSI2019686} and is widely used. Despite their advantages, these methods lack the controllable convergence properties of traditional numerical discretizations, and are biased towards the lower frequency features of the solutions \cite{wang2022and,rahaman2019spectral,krishnapriyan2021characterizing}.

Hybridization frameworks seek to combine the performance of neural network inference on modern accelerated hardware with the guaranteed accuracy of traditional discretizations developed by the scientific community. The hybridization efforts are algorithmic or architectural.

One important algorithmic method is the deep Galerkin method (DGM) \cite{SIRIGNANO20181339}, a neural network extension of the mesh-free Galerkin method where the solution is represented as a deep neural network rather than a linear combination of basis functions. Being mesh-free, it enables the solution of high-dimensional problems by training the neural network model to satisfy the differential operator and its initial and boundary conditions on a randomly sampled set of points, rather than on an exponentially large grid. Although the number of points can be very large in higher dimensions, the training is done sequentially on smaller batches of data points and second-order derivatives are calculated by a scalable Monte Carlo method. Another important algorithmic method is the deep Ritz \cite{yu2018deep}. It implements a deep neural network approximation of the trial function that is constrained by numerical quadrature rules for the variational functional, followed by stochastic gradient descent.

Architectural hybridization is based on differentiable numerical linear algebra. One emerging class involves implementing differentiable finite discretization solvers and embedding them in the neural architectures that enable application of end-to-end differentiable gradient based optimization. Differentiable solvers have been developed in JAX \cite{jax2018github} for fluid dynamic problems, e.g. \texttt{Phi-Flow} \cite{holl2020phiflow}, \texttt{JAX-CFD} \cite{Kochkov2021-ML-CFD}, and \texttt{JAX-FLUIDS} \cite{bezgin2022jax}. These methods are suitable for inverse problems where an unknown field is modeled by the neural network, while the model influence is propagated by the differentiable solver into a measurable residual \cite{pakravan2021solving, dal2020data, lu2020extracting}. We also note the classic strategy for solving inverse problems, the adjoint method, to obtain the gradient of the loss without differentiation across the solver \cite{berg2017neural}; however, deriving analytic expression for the adjoint equations can be tedious or impractical. Other notable use of differentiable solvers is to model and correct for the solution errors of finite discretizations \cite{um2020solver}, and learning and controlling differentiable systems \cite{de2018end, holl2020learning}.

Neural networks are not only universal approximators of continuous functions, but also of nonlinear operators \cite{chen1995universal}. Although this fact has been leveraged using data-driven strategies for learning differential operators by many authors \cite{lu2019deeponet,bhattacharya2020model,li2020neural,li2020fourier}, researchers have demonstrated the ability of differentiable solvers to effectively train nonlinear operators without any data in a completely physics-driven fashion, see section on learning inverse transforms in \cite{pakravan2021solving}.

We propose a novel framework for solving PDEs based on deep neural networks, that enables lifting any existing mesh-based finite discretization method off of its underlying grid and extend it into a mesh-free and embarrassingly parallel method that can be applied to high dimensional problems on unstructured random points. In addition, discontinuous solutions can be readily considered.

\section{Problem statement}
We illustrate our approach by considering a closed irregular interface ($\rm \Gamma$) that partitions an interior ($\rm \Omega^-$) and an exterior ($\rm \Omega^+$) subdomain (see figure \ref{fig:nbm}). The coupled solution $\rm u^\pm\in \Omega^\pm$ satisfy the Helmholtz equation $k^{\pm}u^{\pm} - \nabla \cdot (\mu^{\pm}\nabla u^\pm)=f^{\pm}$ with jump conditions $[u]=\alpha$ and $[\mu \partial_{\mathbf{n}}u]=\beta$,
where $f^\pm$, $\mu^\pm$ and $k^\pm$ are spatially varying coefficients. For simplicity, we consider Dirichlet conditions at the boundary of a cubic computational domain.

This class of problems describes core components of diffusion dominated processes, where sharp irregular interfaces regulate transport across regions with different properties. Examples include Poisson-Boltzmann equation for describing electrostatic properties of biomolecules with jump in dielectric permittivities \cite{sharp1990calculating,MirzadehPB}, in electroporation of cell aggregates with nonlinear membrane jump conditions \cite{mistani2019parallel}, or in epitaxial growth \cite{MISTANI2018150}. Other important applications are found in solidification of multicomponent alloys in additive manufacturing \citep{theillard2015sharp,bochkov2021sharp}, directed self-assembly of diblock copolymers for next generation lithography \cite{galatsis2010patterning,ouaknin2018level,bochkov2021non}, and multiphase flows with phase change.

\section{Scalable and Mesh-Freeing Neuro-Symbolic Differential Solver}

\begin{figure}
\centering
\includegraphics[width=\textwidth]{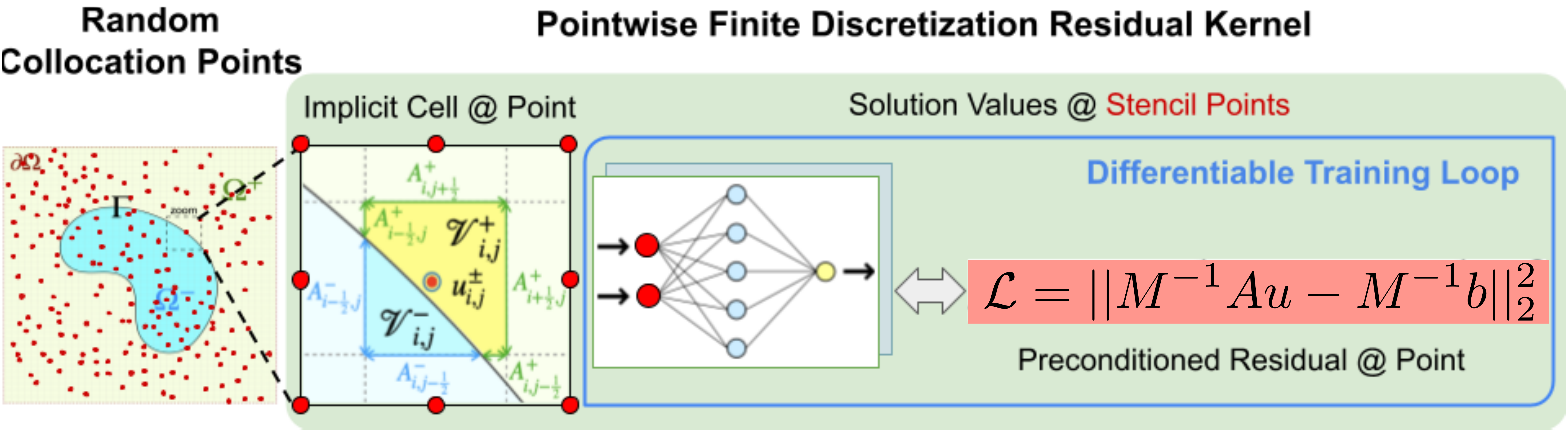}
\caption{Neural Bootstrapping Method (NBM). NBM kernels compute the residual contribution at each collocation point per GPU thread. Kernel operations involve placing implicit cells at multi-resolutions and calculating interface-cell crossings. Point-wise residuals are locally preconditioned ($M^{-1}$ matrix). $A,b$ are the left/right-hand sides of the linear system obtained by finite discretization method. $u$ is the vector of solution values predicted by the neural network model at the stencil points.}
\label{fig:nbm}
\end{figure}
We use neural networks as surrogates for the solution function that are iteratively adjusted to minimize discretization residuals at a set of randomly sampled points and at arbitrary resolutions. The key idea is that neural networks can be evaluated over vertices of any discretization stencils centered at any point in the domain, effectively emulating any structured mesh without actually materializing it. Therefore, we use neural networks to create mesh-free neural differential solvers from existing mesh-based discretization methods . We call this the Neural Bootstrapping Method (NBM) (figure \ref{fig:nbm}).

Numerical methods offer guaranteed accuracy and controllable convergence properties for the training of neural network surrogate models. NBM offers a straightforward path for applying mesh-based numerical methods on random points. This is an important capability for augmenting observational data in the training pipelines. NBM is also a highly parallelizable strategy and the point-wise nature of its kernels is ideally suited for GPU-accelerated computing. The difficult problem of multi-GPU parallel solutions of differential systems is thus reduced to the much simpler task of data-parallel training using existing machine learning frameworks. Data parallelism involves distributing collocation points across multiple processors to compute gradient updates and then aggregating these locally computed updates \cite{shallue2018measuring}.

Here we bootstrap the numerical scheme proposed by \cite{BOCHKOV2020109269} to solve the Helmholtz equation, hence the loss function is 
with the coefficients described in figure \ref{fig:nbm}. The operations in differentiable NBM kernels are strictly local. NBM starts by placing implicit compute cells of a specified resolutions at the collocation point. At discontinuities, two neural networks are used to represent the solution on each side of $\Gamma$ and a separate coarse mesh encapsulates a level-set function \citep{osher1988fronts} that provides the necessary geometric information for the numerical kernel and the preconditioner (denoted $M^{-1}$). The numerical kernel is applied on the compute cell where the solution values are evaluated by the neural network. Each kernel contributes a local $L^2$-norm residual $r_p = ||M^{-1}A u - M^{-1}b||$ at one point $p$. Preconditioning balances the relative magnitude of contributions from all points before aggregating the residuals to form a global loss value. Finally, gradient based optimization methods used in machine learning, \textit{e.g.} \texttt{Adam} optimizer \cite{kingma2014adam}, are applied to adjust neural network parameters. The automatic differentiation loop passes across the NBM kernels.

\section{Numerical results}
\subsection{Convergence and accuracy}

We consider a sphere centered at the origin with radius 0.5 in a domain $[-1, 1]^3$, a discontinuous exact solution of $u^-=e^{z}$ and $u^+=\cos(x)\sin(y)$. In addition to the jump in solution, we consider a jump in the variable diffusion coefficient to be $\mu^-=y^2 \ln(x+2) + 4$ and $\mu^+=e^{-z}$. Table \ref{tab:convergence} reports convergence in the $L^\infty$-norm. Two neural networks were used to represent solutions inside and outside the sphere with 5 hidden layers and 10 sine-activated neurons each, for a total of $982$ trainable parameters.

\begin{table}[ht]
\begin{center}
\caption{Convergence and overall time to solution for our JAX implementation, with $10,000$ epochs in each case, and initial compilation time. Measurements are on a single NVIDIA A6000 GPU. } \label{tab:convergence}

\begin{tabular}{ccccccc}
\toprule
 & \multicolumn{2}{c}{RMSE}& \multicolumn{2}{c}{$L^\infty$} & \multicolumn{2}{c}{GPU Statistics} \\
\cmidrule(r){2-3}\cmidrule(r){4-5}\cmidrule(r){6-7}
$\rm N_{x,y,z}$   &   Solution    &   Order   &   Solution   &   Order & t (sec/epoch) & VRAM (GB)\\
\midrule
$2^3$ & $3.7\times 10^{-2}$ &  -        & $3.25\times 10^{-1}$  &   -     & $0.0306$ & $1.05$ \\
$2^4$ & $7.1\times 10^{-3}$ &  $2.38$   & $1.10\times 10^{-1}$   & $1.56$  & $0.056$  & $1.72$ \\ 
$2^5$ & $5.9\times 10^{-3}$ &  $0.27$   & $8.36 \times 10^{-2}$ & $0.4$   & $0.053$  & $2.15$ \\ 
$2^6$ & $4.1\times 10^{-3}$ &  $0.53$   & $6.44\times 10^{-2}$  & $0.38$  & $0.287$  & $5.57$ \\ 
$2^7$ & $2.64\times 10^{-3}$&  $0.64$   & $3.53\times 10^{-2}$  & $0.87$  & $2.125$ & $32.1$ \\ \bottomrule
\end{tabular}
\end{center}
\end{table}

\subsection{Time complexity and parallel scaling on GPU clusters}\label{sec:scaling}

We simulate a Helmholtz problem with discontinuities on the Dragon geometry presented in \cite{curless1996volumetric}. We consider a solution with jumps across the dragon's surface and a spatially varying diffusion coefficient. The results are shown in figure \ref{fig:dragon}, with a $\rm L^\infty$-error of $0.5$ and RMSE of $0.06$ after 1000 epochs on a base resolution of $64^3$ and implicitly refined onto multi-resolutions $128^3,256^3,512^3$. The neural network pair have only 1 hidden layer with 100 sine-activated neurons, although investigating more complex networks (transformers, symmetry preserving, \textit{etc}.) would likely improve accuracy.

In table \ref{tab:scaling} we report scaling results on NVIDIA A100 GPUs at four base resolutions with three levels of implicit refinement. We used a batchsize of $32\times 32\times 16$ in all cases. At fixed number of GPUs, training time scales linearly (\textit{i.e.}, optimal scaling) with the number of grid points. At a fixed resolution, increasing the number of GPUs accelerates training roughly with $\rm epoch\ time \sim 1/\sqrt{\text{\# }GPUs}$, although the advantage is more effective at higher resolutions. Compile time increases with resolution and decreases with number of GPUs. A maximum grid size of $1024^3$ at multi-resolutions $1024^3,\ 2048^3, 4096^3,\ 8192^3$ was simulated on one NVIDIA DGX with 8 A100 GPUs. 

\begin{table}[ht]
\begin{center}
\caption{Scaling test. Time per epoch (sec) and JAX compile time for different configurations.} \label{tab:scaling}
\begin{tabular}{ccccccccc}
\toprule
base resolution: & \multicolumn{2}{c}{$64^3$}& \multicolumn{2}{c}{$128^3$} & \multicolumn{2}{c}{$256^3$}  & \multicolumn{2}{c}{$512^3$}\\
\cmidrule(r){1-1}\cmidrule(r){2-3}\cmidrule(r){4-5}\cmidrule(r){6-7}\cmidrule(r){8-9}
$\rm A100\ GPUs$ &  epoch  &  compile   &  epoch   & compile &    epoch & compile  &  epoch    &   compile \\
\cmidrule(r){1-1}\cmidrule(r){2-3}\cmidrule(r){4-5}\cmidrule(r){6-7}\cmidrule(r){8-9} 
$1$        & $0.908$ &  $9.027$   & $6.960$  & $9.288$ & $55.287$ & $12.164$ &  $438.45$  &  $49.020$ \\
$2$        & $0.657$ &  $7.575$   & $5.893$  & $7.823$ & $47.360$ & $10.045$ &  $378.98$  &  $39.815$ \\ 
$4$        & $0.405$ &  $7.480$   & $3.629$  & $7.863$ & $28.261$ & $9.129$  &  $226.73$  &  $27.405$ \\ 
$8$        & $0.384$ &  $7.983$   & $3.340$  & $7.901$ & $26.799$ & $9.154$  &  $204.88$  &  $20.632$\\ \bottomrule
\end{tabular}
\end{center}
\end{table}

\section{Conclusion}
We presented a neural bootstrapping method and applied it to the problem of solving elliptic PDEs with discontinuities. NBM is a differentiable computing method that creates scalable and mesh-free numerical methods from mesh-based finite discretizations. It represents the solution by training neural networks using automatic differentiation of the discretized residual at collocation points. We implemented the method using JAX and showed accuracy and parallel scaling in three dimensions.

\begin{ack}
This work has been partially funded by ONR N00014-11-1-0027.
\end{ack}

\section*{Broader Impact}

The work presented provides a systematic framework for creating neural-based solvers for partial differential equations with potential jump conditions across sharp interfaces, which describe a plethora of important systems in the physical and life sciences. While the driving motivation for our work is the development of algorithms for multi-scale simulations to accelerate drug discovery and formulation pipelines in the pharmaceutical industry, we believe that the methodology will positively impact other industries that face similar computational challenges. We do not believe that the methodology we have introduced will have an adverse effect to society or will have negative ethical implications.



\bibliographystyle{abbrv}
\addcontentsline{toc}{section}{\refname}
\bibliography{references}

\section*{Checklist}

\begin{enumerate}

\item For all authors...
\begin{enumerate}
  \item Do the main claims made in the abstract and introduction accurately reflect the paper's contributions and scope?
    \answerYes{}
  \item Did you describe the limitations of your work?
    \answerYes{see section \ref{sec:scaling} on exploring more complex neural network architectures for improving accuracy.}
  \item Did you discuss any potential negative societal impacts of your work?
    \answerYes{See section on Broder Impact.}
  \item Have you read the ethics review guidelines and ensured that your paper conforms to them?
    \answerYes{}
\end{enumerate}

\item If you are including theoretical results...
\begin{enumerate}
  \item Did you state the full set of assumptions of all theoretical results?
    \answerNA{}
        \item Did you include complete proofs of all theoretical results?
    \answerNA{}
\end{enumerate}

\item If you ran experiments...
\begin{enumerate}
  \item Did you include the code, data, and instructions needed to reproduce the main experimental results (either in the supplemental material or as a URL)?
    \answerNo{The code is proprietary.}
  \item Did you specify all the training details (e.g., data splits, hyperparameters, how they were chosen)?
    \answerYes{}
        \item Did you report error bars (e.g., with respect to the random seed after running experiments multiple times)?
    \answerNo{}
        \item Did you include the total amount of compute and the type of resources used (e.g., type of GPUs, internal cluster, or cloud provider)?
    \answerYes{}
\end{enumerate}

\item If you are using existing assets (e.g., code, data, models) or curating/releasing new assets...
\begin{enumerate}
  \item If your work uses existing assets, did you cite the creators?
    \answerNA{}
  \item Did you mention the license of the assets?
    \answerNA{}
  \item Did you include any new assets either in the supplemental material or as a URL?
    \answerNA{}
  \item Did you discuss whether and how consent was obtained from people whose data you're using/curating?
    \answerNA{}
  \item Did you discuss whether the data you are using/curating contains personally identifiable information or offensive content?
    \answerNA{It is a completely physics driven work without any data.}
\end{enumerate}

\item If you used crowdsourcing or conducted research with human subjects...
\begin{enumerate}
  \item Did you include the full text of instructions given to participants and screenshots, if applicable?
    \answerNA{}
  \item Did you describe any potential participant risks, with links to Institutional Review Board (IRB) approvals, if applicable?
    \answerNA{}
  \item Did you include the estimated hourly wage paid to participants and the total amount spent on participant compensation?
    \answerNA{}
\end{enumerate}

\end{enumerate}

\end{document}